\documentclass[12pt]{amsart}
\usepackage{latexsym,enumerate}
\usepackage{amsmath,amsthm,amsopn,amstext,amscd,amsfonts,amssymb,mathrsfs}
\usepackage{color}
\usepackage{graphicx}

\numberwithin{equation}{section}

\newtheorem{theorem}{Theorem}
\newtheorem{lemma}{Lemma}
\newtheorem{corollary}{Corollary}
\newtheorem{proposition}{Proposition}

\numberwithin{theorem}{section}
\numberwithin{corollary}{section}
\numberwithin{lemma}{section}
\numberwithin{definition}{section}
\numberwithin{proposition}{section}
\numberwithin{remark}{section}

\newcommand{\medint}{-\kern  -,375cm\int}
\newcommand{\dint}{\displaystyle\int}

\begin{document}

\title[]{On  weighted isoperimetric inequalities \\with non-radial densities }

\author[A. Alvino]{A. Alvino$^1$}
\author[F. Brock]{F. Brock$^2$}
\author[F. Chiacchio]{F. Chiacchio$^1$}
\author[A. Mercaldo]{A. Mercaldo$^1$}
\author[M.R. Posteraro]{M.R. Posteraro$^1$}


\setcounter{footnote}{1}
\footnotetext{Universit\`a di Napoli Federico II, Dipartimento di Matematica e Applicazioni ``R. Caccioppoli'',
Complesso Monte S. Angelo, via Cintia, 80126 Napoli, Italy;\\
e-mail: {\tt angelo.alvino@unina.it, fchiacch@unina.it,  mercaldo@unina.it, posterar@unina.it}}

\setcounter{footnote}{2}
\footnotetext{South China University of Technology, International School of 
Advanced Materials, Wushan Campus 381, Wushan Road, 
Tianhe District, Guangzhou, P.R. China, 510641, email: fbrock@scut.edu.cn
\medskip

The paper was partially supported by the grants GNAMPA
}

\begin{abstract} 
We consider a class of isoperimetric problems on $\mathbb{R}^{N}_{+} $ where
the volume and the area element carry two different weights of the type $|x|^lx_N^\alpha$. We solve them in a special case while a more detailed study is contained in \cite{ABCMP2}.  Our results  imply a weighted Polya-Sz\"ego principle and a priori estimates for weak solutions to a class of boundary value problems for degenerate elliptic equations\medskip

\rm 
\\[0.1cm]
\rm 
\end{abstract}
\maketitle

\section{Introduction }
The study of isoperimetric problems when volume and perimeter carry two different weights 
has  recently attracted the attention of many authors. 
A remarkable example is obtained when
 the two weights are  two
different powers of the distance from the origin.
More precisely given two real numbers $k$ and $l$ the problem is to find the set $G$ in $\mathbb R^{N}$
which minimizes the weighted perimeter $\displaystyle\int_{\partial G } |x|^k  \, {\mathcal H}_{N-1} (dx) $ 
once it is prescribed the weighted measure $\displaystyle\int_{ G } |x|^l  \, dx. $ 
The interest for such a problem is due to the fact that its solution allows to get, for instance,
 the knowledge of the best constants in the well-known Caffarelli-Kohn-Nirenberg inequalities as well as 
the radiality of the corresponding minimizers.
Several partial results have been obtained in this framework (see, e.g., 
 \cite{ABCMP}, \cite{BBMP2},  \cite{C},   \cite{diGiosia_etal},  \cite{DHHT}, \cite{Howe},  \cite{Mo}   and the references therein)) and a
complete solution is contained in the recent paper  \cite{diGiosia_etal} where the authors reduce the problem
into a two-dimensional one by means of spherical symmetrization.
\\
The problem we address here is the following:

Given  $k,l \in \mathbb{R}$, $\alpha > 0$,
\medskip
\noindent {\sl Minimize $\displaystyle\int_{\partial \Omega } |x|^k x_N^\alpha \, {\mathcal H}_{N-1} (dx) $ 
among all smooth sets 
$\Omega \subset \mathbb{R} ^{N}_{+}$ satisfying $\dint_{\Omega } |x|^lx_N^\alpha \, dx =1$}
\medskip
where $\mathbb{R}^{N} _{+}  := \{ x \in \mathbb{R}^N :\, x_N >0\} $.  
\noindent In the present paper we give a partial result, while a more extensive study is contained in \cite{ABCMP2}.

\noindent Let $B_R$ denote the ball of $\mathbb R^N$ of radius $R$ centered at the origin
and  let $B$ and $\Gamma$ denote the Beta and the Gamma function, respectively.

\noindent The main  result of this paper, proved in Section 4,  is the following
\begin{theorem}
\label{maintheorem}
Let $N\in \mathbb{N} $, $N\geq 2$, 
$k,l \in \mathbb{R} $, $\alpha > 0$ and $l+N+\alpha >0$. 
Further assume that 
\begin{equation*}\label{i}
l+1\leq k .
\end{equation*}
Then  
\begin{equation}
\label{mainineq}
\dint_{\partial \Omega } |x|^k x_N^\alpha\, {\mathcal H}_{N-1} (dx)
\geq 
C_{k,l,N, \alpha} ^{rad}  
\left( 
\dint_{\Omega } |x|^lx_N^\alpha\, dx 
\right) 
^{(k+N+\alpha-1)/(l+N+\alpha) } , 
\end{equation}
for all smooth sets $\Omega $ in 
$\mathbb{R}^N _+  $,
where 
\begin{eqnarray*}
\label{defCkl}
C_{k,l,N, \alpha} ^{rad} & := &
\frac{\dint_{\partial B_1 } |x|^k x_N^\alpha\, {\mathcal H}_{N-1} (dx)}
{\left( \dint_{B_1 \cap  \mathbb{R}^N _+}
 |x|^l x_N^\alpha\, dx \right ) ^{(k+N+\alpha-1)/(l+N+\alpha) } } 
\\
 &= &
\nonumber
\left( l+\alpha +N\right) ^{\frac{k+N+\alpha -1}{l+N+\alpha }}\left(
B\left( \frac{N-1}{2},\frac{\alpha +1}{2}\right) \frac{\pi ^{\frac{N-1}{2}}}{
\Gamma \left( \frac{N-1}{2}\right) }\right) ^{\frac{l-k+1}{l+N+\alpha }}.
\end{eqnarray*}
Equality in (\ref{mainineq}) holds  if $\Omega =B_R\cap\mathbb{R}_+^N$.
\end{theorem}

\noindent

Note that the weights we consider are not radial and therefore it seems not trivial to use
 spherical symmetrization. For this reason our approach is different from the one used in  \cite{diGiosia_etal}.
 
In \cite{ABCMP2} also the cases where 
$k\leq l+1$ and $ l\frac{N+\alpha-1}{N+\alpha} \leq k\leq 0$; 
or 
$N\geq 2$, $ 0\leq k\leq l+1$ and 
\begin{equation*}\label{l_1N3}
l\le l_1 (k,N,\alpha )  :=  \frac{(k+N+\alpha-1)^3 }{(k+N+\alpha-1)^2 - \frac{(N+\alpha-1)^2 }{N+\alpha} } -N -\alpha\,.
\end{equation*}
are faced and, depending on the regions where the three parameters lie, we use different methods.
In \cite{ABCMP2} we also provide some necessary conditions on $k$, $l$ and $\alpha$  such that 
the half-ball centered at the origin is an isoperimetric set. 

Several applications of our isoperimetric problems are proved in \cite{ABCMP2} such as, for example, a Polya-Sz\"ego type principle or  the best constants in Caffarelli-Kohn-Nirenberg type inequalities.
In this note we just prove a special case of a Polya-Sz\"ego type principle from which we deduce a Poincar\'e type inequality. Such inequality allows to ensure, via Lax-Milgram Theorem, the existence of a weak solution to a boundary value problem for a class of degenerate elliptic equations. For weak solutions to such problems we state a comparison result which can be proved by using the same method of Talenti \cite{T} or Alvino and Trombetti  \cite{AT}.

\section{Notation and preliminary results  } 

Throughout this article  $N$ will denote a natural number with $N\geq 2$,  
$k$ and $l $ are real numbers, while $\alpha$ is a positive number
such that  
\begin{equation}
\label{ass1}
 l+N+\alpha>0 .
\end{equation}

\noindent 
As usual we denote
\begin{eqnarray*}
\mathbb{R}^{N}_{+}
& := & 
\left\{ x \in \mathbb{R}^{N}: x_N  >0 \right\}, 
\end{eqnarray*}
and for any $R>0$
\begin{eqnarray*}
B_{R} 
& := & 
\left\{ x\in \mathbb{R}^{N}:\left\vert x\right\vert <R\right\} \\
B_{R}^{+}  
& := & 
B_{R}  \cap \mathbb{R}^{N}_{+}.
\\
\end{eqnarray*}
Furthermore, 
we will set
\begin{equation*}
\omega _N  :=  {\mathcal L}^N (B_1 )
\end{equation*}
and
\begin{equation*}
\kappa (N, \alpha )  := \int_{\partial B_1^{+} } x _{N}^{\alpha }
{\mathcal H} _{N-1}(dx)\,
\end{equation*}
where ${\mathcal L} ^N $  denotes  the $N$-dimensional Lebesgue measure and 
 ${\mathcal H} _{N-1} $  the  $(N-1)$-dimensional Hausdorff measure.

\noindent In the Appendix of \cite{BCM3} the value of $\kappa (N, \alpha ) $ is computed and it is proved that
\begin{equation*}
\label{measSN-1+}
\kappa (N, \alpha ) = B\left( \frac{
N-1}{2},\frac{\alpha +1}{2}\right) 
\frac{\pi ^{\frac{N-1}{2}}}{\Gamma \left( 
\frac{N-1}{2}\right) },
\end{equation*}
where $B$ and $\Gamma$ are the Beta function and the Gamma function,
respectively.

\noindent Next we define a measure $\mu _{l, \alpha}$ by 
\begin{equation*}
d\mu _{l, \alpha}(x)=|x|^{l} x_N^\alpha\,dx.  
\label{dmu}
\end{equation*}
If $M \subset $ ${\mathbb R}^{N}_{+}$ is a  measurable set with finite 
$\mu _{l, \alpha} $-measure, then we define $M^{\star }$, the 
\\
$\mu_{l,\alpha }$-symmetrization of $M$,
 as follows
\begin{equation}
M^{\star } := B_{R^\star}^{+} 
\hspace{.2 cm}
 \text{with }
 R^\star:
\mu_{l, \alpha} 
\left( B_{R^\star}^{+}    \right) =
\mu _{l, \alpha} \left( M\right) = \int_M d\mu _{l, \alpha} (x) .  
\label{mu_(M)}
\end{equation}
A straightforward calculation gives 
$$
R^\star=\left (\frac{\mu_{l,\alpha}(M)}{\mu_{l,\alpha}(B_1^+)} \right )^{\frac{1}{l+\alpha+N}}
$$
where 
$$
\mu_{l,\alpha}(B_1^+)=\frac{\kappa (N, \alpha )}{l+\alpha+N}\,.
$$
If $u: M \rightarrow \mathbb{R} $ is a measurable function such that
$$
\mu_{l, \alpha} \left( \left\{ |u(x)|>t\right\} \right) <\infty \qquad \forall t>0,
$$
then  $u^*$ denotes the weighted decreasing rearrangement of $u$,  which is given by
\begin{equation*}
u^*(s)=\sup \left\{ t\geq 0:\mu_{l, \alpha}
 \left( \left\{ |u(x)| >t\right\} \right)
>s
 \right\}\,, \qquad\qquad s\in ]0,\mu_{l, \alpha}(M)] \,,  
\label{u_star}
\end{equation*}
and $u^{ \star }$ denotes the weighted Schwarz symmetrization of $u$, or
in short, the \\
$\mu_{l, \alpha} -$symmetrization of $u$, which is given by
$$
u^{ \star }(x)= u^*(\mu_{l,\alpha}(B^+_{|x|}))= u^*(\mu_{l,\alpha}(B^+_1)|x|^{l+\alpha+N})\,, \qquad x\in M^\star .
$$ 
\noindent Note that $u^{\star }$ is a radial and radially non-increasing function.

\noindent Moreover if $M$ is a measurable set in $\mathbb{R}^N _{+}$ with finite $\mu _{l,\alpha} $-measure, then 
$$
\left( \chi _M \right) ^{\star} = \chi _{M^{ \star }} \,.
$$
Here and in the following for any   measurable set $M$ in $\mathbb{R}^N _{+}$,  $\chi _M $  denotes its characteristic function.

The {\sl $\mu _{k, \alpha}$--perimeter\/} of a measurable set $M $ is given by 
\begin{equation*}
P_{\mu _{k, \alpha}}(M ):=\sup \left\{ \int_{M }\mbox{div}\,\left(x_N^\alpha |x|^{k}
\mathbf{v}\right) \,dx:\,\mathbf{v}\in C_{0}^{1}(\mathbb{R}^N ,\mathbb{R}^{N}),\,|
\mathbf{v}|\leq 1\mbox{ in }\, M \right\} .
\end{equation*}

\noindent It is well-known that the above 
\textsl{distributional definition} of
weighted perimeter is equivalent to the following

\begin{equation*}
P_{\mu_{k} }(M )
= \left\{ 
\begin{array}{ccc}
\displaystyle\int_{\partial \Omega }|x|^{k} \, {\mathcal H} _{N-1}(dx)
 & \mbox{ if } & 
\partial \Omega  \mbox{ is } (N-1)-\mbox{rectifiable } \\ 
&  &  \\ 
+ \infty \qquad
 & \mbox{ otherwise.} & 
\end{array}
\right.
\end{equation*}


Let $\Omega \subset \mathbb{R} ^{N}_{+}$ and $p\in \left[ 1,+\infty \right) $.
We will denote by $L^{p}(\Omega ,d\mu
_{l, \alpha})$ the space of all Lebesgue measurable real valued functions $u$ such
that
\begin{equation*}
\left\Vert u \right\Vert 
_{L^{p}(\Omega ,d\mu _{l, \alpha})}
:=\left( 
\int_{\Omega
}
\left\vert 
u\right\vert 
^{p}d\mu _{l, \alpha} (x) 
\right) ^{1/p}
<+\infty .
\label{Norm_Lp}
\end{equation*}
\\  
By $W^{1,p}(\Omega ,d\mu _{l, \alpha})$ we denote the weighted Sobolev space
consisting of all functions which together with their weak derivatives $u_{x_{i}}$, ($i=1,...,N$), 
belong to $L^{p}(\Omega ,d\mu _{l, \alpha})$.
This space will be equipped with the norm
\begin{equation*}
\left\Vert u\right\Vert _{W^{1,p}(\Omega ,d\mu _{l, \alpha})}:=\left\Vert
u\right\Vert _{L^{p}(\Omega ,d\mu _{l, \alpha})}+\left\Vert \nabla u\right\Vert
_{L^{p}(\Omega ,d\mu _{l, \alpha})}.  
\label{Norm_Wp}
\end{equation*}
Let $X$ be the set of all functions $u\in C^{1}(\bar \Omega)$ that vanish in a neighborhood of $\partial \Omega \setminus \{ x_N=0\}$. Then let $V^2(\Omega, d\mu_{l,\alpha})$ be the closure of $X$ in the norm of $W^{1,2}(\Omega,d\mu_{l,\alpha} )$.

\section{The isoperimetric problem }

In order to formulate the isoperimetric problem  we introduce two functionals ${\mathcal R}_{k,l,N,\alpha}$ and ${\mathcal Q}_{k,l,N,\alpha}$.

\noindent 
If $M $ is any measurable subset of $\mathbb R^{N}_{+}$, with $0<\mu _{l,\alpha} (M)<+\infty $, we set
\begin{equation*}
\label{rayl1}
{\mathcal R}_{k,l,N, \alpha} (M) := 
\frac {P_{ \mu_{k , \alpha}}  (M) }
{  \left( \mu_{l,\alpha} (M) \right)^{(k+N+\alpha-1)/(l+N+\alpha)} }. 
\end{equation*}
Note that
\begin{equation*}
\label{Rklsmooth}
{\mathcal R} _{k,l,N,\alpha} (M ) = 
\frac{
\dint_{\partial M }x_N^\alpha |x|^k \, {\mathcal H}_{N-1} (dx) 
}{
\left( \dint_{M }x_N^\alpha |x|^l \, dx \right) ^{(k+N+\alpha-1)/(l+N+\alpha)} } 
\end{equation*}
if the set $M$ is smooth.

If $u\in C_0 ^1  (\mathbb{R} ^N_+ )\setminus \{ 0\} $, we set
\begin{equation*}
\label{rayl2}
{\mathcal Q}_{k,l,N, \alpha} (u ) := \frac{\dint_{\mathbb{R} ^N_+ }x_N^\alpha |x|^k |\nabla u| \, dx}{ 
\left( \dint_{\mathbb{R} ^N_+ } x_N^\alpha|x|^l |u| ^{(l+N+\alpha)/(k+N+\alpha-1)} \, dx \right) ^{(k+N+\alpha-1)/(l+N+\alpha)}}. 
\end{equation*}

Finally, we define 
\begin{equation*}
\label{isopco}
C_{k,l,N, \alpha}^{rad} := {\mathcal R}_{k,l,N, \alpha}(B_1 \cap {\mathbb{R} ^N_+ }).
\end{equation*} 
We study the following isoperimetric problem: 
\\[0.5cm]
{\sl Find the constant $C_{k,l,N, \alpha} \in [0, + \infty )$, such that}
\begin{equation*}
\label{isopproblem}
C _{k,l,N, \alpha} := \inf \{ {\mathcal R}_{k,l,N, \alpha} (M):\, 
\mbox{{\sl $M$ is measurable with $0<\mu _{l,\alpha} (M) <+\infty $}} \}. 
\end{equation*}   
Moreover, we are interested in conditions on $k$, $l$ and $\alpha$ such that
\begin{equation}
\label{isoradial}
{\mathcal R}_{k,l,N, \alpha} (M) \geq {\mathcal R}_{k,l,N, \alpha} (M^{ \star} ) 
\end{equation}
holds for all measurable sets $M\subset {\mathbb{R} ^N_+ }$ with  $ 0<\mu _{l,\alpha}(M)<+\infty $. 
\\[0.1cm]
Let us begin with some standard remarks.
\\
If $M$ is a measurable subset of $\mathbb R^{N}_{+}$ with finite $\mu _{l,\alpha}$-measure and $\mu_{k,\alpha} $-perimeter,
 then there exists a sequence of smooth sets 
$\{ M_n \} $ such that $\displaystyle\lim_{n\to \infty } \mu _{l,\alpha} (M_n \Delta M) =0$ and 
$\displaystyle\lim_{n\to \infty } P_{\mu _{k,\alpha} } (M_n ) = P_{\mu _{k,\alpha}} (M)$. 
This property is well-known for Lebesgue measure (see for instance 
\cite{G}, Theorem 1.24) 
and its proof carries over to the weighted case. This implies that we also have 
\begin{equation}
\label{CklNsmooth}
C_{k,l,N, \alpha} = \inf \{ {\mathcal R}_{k,l,N, \alpha} (\Omega ):\, \Omega \subset \mathbb{R} ^{N}_+, \, \Omega  
\mbox{ smooth} \} .
\end{equation}  

\noindent Actually in \cite{ABCMP2} it is proved that for the value of $k,l$ such that $k\le l+1$, the constant $C_{k,l,N, \alpha}$ is related to the functional ${\mathcal Q}_{k,l,N, \alpha }$ as follows
$$
C_{k,l,N, \alpha} = \inf \{ {\mathcal Q}_{k,l,N, \alpha} (\Omega ):\, u\in C_0 ^1 (\mathbb{R}^N_+ )\setminus \{ 0\}  \} .
$$
The functionals ${\mathcal R}_{k,l,N, \alpha } $ and ${\mathcal Q}_{k,l,N,\alpha } $ 
have the following homogeneity properties, 
\begin{eqnarray*}
\label{hom1}
 {\mathcal R}_{k,l,N, \alpha } (M ) & = & {\mathcal R}_{k,l,N, \alpha } (tM ) ,
\\
{\mathcal Q}_{k,l,N, \alpha } (u) & = & {\mathcal Q}_{k,l,N,\alpha } (u^t ),
\end{eqnarray*}
where  $t>0$, $M $ is a measurable set with $0<\mu_{l, \alpha} (M)<+\infty $, 
$u\in C_0 ^1 (\mathbb{R}^N_+ )\setminus \{ 0\}$,  
\\
$tM := \{tx:\, x\in M \} $ 
and $u^t (x):= u(tx) $, ($x\in \mathbb{R} ^N_+ $), and there holds   
\begin{equation*}
\label{isopconst2}
C_{k,l,N, \alpha} ^{rad} = {\mathcal R}_{k,l,N, \alpha} (B_1^{+} ).
\end{equation*}  
Hence we have that 
\begin{equation}
\label{relCC}
C_{k,l,N, \alpha} \leq C_{k,l,N, \alpha} ^{rad} ,
\end{equation}  
and (\ref{isoradial}) holds if and only if 
\begin{equation}\label{isop1}
C_{k,l,N,\alpha } = C_{k,l,N,\alpha } ^{rad} .
\end{equation}  
Finally, we recall the following weighted isoperimetric inequality proved, for example, in \cite{BCM2} 
(see also \cite{BCM3, XR, XRS, MadernaSalsa})

\begin{proposition}\label{BCM2} 
For all measurable sets $M\subset \mathbb{R} ^N_+$, with $0< \mu _{0, \alpha}  (M)<+\infty $, the following inequality holds true
\begin{equation*}\label{isopclass}
{\mathcal R} _{0,0,N, \alpha} (M) := 
\frac {P_{ \mu_{0, \alpha}}  (M) }
{  \left( \mu_{0,\alpha} (M) \right)^{(N+\alpha-1)/(N+\alpha)} } \geq C_{0,0,N, \alpha} ^{rad}:= 
\frac {P_{ \mu_{0, \alpha}}  (M^{ \star}) }
{  \left( \mu_{0,\alpha} (M^{ \star}) \right)^{(N+\alpha-1)/(N+\alpha)} } \,,
\end{equation*}
where $M^{\star}=B_{R}^{+} $ with $R$ such that $\mu_{0, \alpha}(M)=\mu_{0, \alpha}(M^{ \star})$
\end{proposition}

We recall that the isoperimetric constant $C_{0,0,N, \alpha} ^{rad}$ is explicitly computed in \cite{BCM2} (see also \cite{MadernaSalsa} for the case $N=2$).

\section{Main result}

This section is devoted to the proof of Theorem \ref{maintheorem} which provides  
 a sufficient condition on $k,l$  and $\alpha$ such that 
$ C_{k,l,N, \alpha} = C_{k,l,N, , \alpha} ^{rad}$ holds, or equivalently,
\begin{equation}
\label{ineqrad}
{\mathcal R}_{k,l,N, \alpha} (M) \geq C_{k,l,N, \alpha}^{rad}  
\quad \mbox{for all measurable sets $M \subset \mathbb R^{N}_{+}$ with $0< \mu _{l, \alpha} (M) <+\infty $.}
\end{equation}
This implies that our main result can be equivalently formulated as follows

\begin{theorem}    
\label{R5}Let $N\in \mathbb{N} $, $N\geq 2$, 
$k,l \in \mathbb{R} $, $\alpha > 0$,  $l+N+\alpha >0$ and  
$l+1\leq k $.
 Then (\ref{isop1}) holds. 

\noindent Moreover, if $l+1 <k$ and
\begin{equation}
\label{M=BR}
{\mathcal R}_{k,l,N, \alpha} (M) = C_{k,l,N, \alpha} ^{rad} \ 
\mbox{ for some measurable set $M \subset \mathbb R^{N}_{+}$ with $0<\mu _l (M)< +\infty $},
\end{equation}  
then $M = B_{R}^{+} $ for some $R>0$. 
\end{theorem}

The proof of this a result relies on Gauss' Divergence Theorem and the following lemmas (see
\cite{KZ}, \cite{BrasPhil} and \cite{ABCMP}).

\begin{lemma}
\label{hardylitt}
Let $l>l' >-N -\alpha$. Then 
\begin{equation} 
\label{hardylitt2}
\frac{\left( \mu _{l, \alpha} (M)  \right) 
^{1/(l+N+\alpha)} 
}{ 
\left( \mu _{l', \alpha} (M)  \right)
^{1/(l'+N+\alpha)}
} 
\geq \frac{\left( \mu _{l, \alpha} (M^{ \star})  \right) 
^{1/(l+N+\alpha)} 
}{ 
\left( \mu _{l', \alpha} (M^{ \star})  \right)
^{1/(l'+N+\alpha)}
} 
\end{equation}
for all measurable sets $M\subset \mathbb{R} ^N_+ $ with  $0<\mu_{l, \alpha}(M)<+\infty $.
Equality holds only for half-balls $B_{R}^{+} $, ($R>0$).
\end{lemma}

{\sl Proof: } Let $M^{ \star} $ be the 
$\mu_{l, \alpha} $-symmetrization of $M$. Then we obtain, using the Hardy-Littlewood inequality,
\begin{eqnarray*}
\mu _{l' , \alpha} (M) =\int_Mx_N^\alpha |x| ^{l'} \, dx & = & \int_{\mathbb{R} ^N_+ }  |x|^{l'-l} \chi _M (x)\, d\mu _{l, \alpha} (x) 
\\
 & \leq & 
\int_{\mathbb{R} ^N _+}  \left( |x|^{l'-l} \right) ^{ \star}  \left( \chi _M  \right) ^{ \star} (x)\, d\mu _{l, \alpha} (x) 
\\
 & = & 
\int_{\mathbb{R} ^N_+ }   |x|^{l'-l}  \chi _{M^{ \star} }  (x)\, d\mu _{l, \alpha} (x) 
\\
 & = & 
\int_{M^{ \star} } x_N^\alpha |x|^{l' }\, dx =\mu _{l', \alpha } (M^{ \star} ).
\end{eqnarray*} 
This implies (\ref{hardylitt2}). 

\noindent Next assume that equality holds in (\ref{hardylitt2}). Then we must have 
$$
\int_M |x|^{l'-l} \, d\mu _{l, \alpha} (x) = \int_{M^{ \star} } |x|^{l'-l} d\mu _{l, \alpha} (x) ,
$$
that is, 
$$
\int_{M\setminus M^{ \star} } |x|^{l'-l} \, d\mu _{l, \alpha} (x) = \int_{M^{ \star} \setminus M} |x|^{l'-l} d\mu _{l, \alpha}  (x) .
$$
Since $l'-l<0$, this means that 
$ \mu _l ( M\Delta M^{ \star} )=0$. The Lemma is proved. 
$\hfill \Box $
\begin{lemma}
\label{rangekl1}
Let $k,l, \alpha$ satisfy (\ref{ass1}). Assume that $l>l' >-N-\alpha$ and 
$C_{k,l,N, \alpha}  = C_{k,l,N, \alpha} ^{rad} $. 
Then we also have 
$C_{k,l',N, \alpha}  = C_{k,l',N, \alpha} ^{rad} $.
Moreover, if  $
{\mathcal R}_{k,l',N, \alpha} (M ) = C_{k,l',N, \alpha} ^{rad} $  
for some measurable set $M\subset \mathbb{R} ^N_+ $, with $0< \mu _{l' , \alpha} (M) <+\infty $,
then $M = B_{R}^{+}$ for some $R>0$.
\end{lemma}
{\sl Proof:} By our assumptions and Lemma \ref{hardylitt}  we have for every measurable set $M$ with 
$0<\mu _{l, \alpha}(M) <+\infty $,  
\begin{eqnarray*}
{\mathcal R}_{k,l',N, \alpha} (M ) & = & {\mathcal R}_{k,l,N, \alpha} (M ) 
\cdot 
\left[
\frac{
\left( 
\mu _{l, \alpha} (M) 
\right) ^{1/(l+N+\alpha)}
}{ 
\left( \mu_{l', \alpha} (M) 
\right) ^{1/(l'+N+\alpha)}
} 
\right] ^{k+N+\alpha-1}
\\
 & \geq & 
C_{k,l',N, \alpha}^{rad},
\end{eqnarray*}
with equality only if 
$M = B^{+}_{R}$ for some $R>0$. 
$\hfill \Box $
\medskip

\noindent{\bf Proof of Theorem \ref{R5}: } We consider two cases. 
\\
Assume $l+1=k$.  Let $\Omega $ be a smooth subset of $\mathbb R^{N}_{+}$ and let  
$\Omega ^{\star}$
be the set defined in \eqref{mu_(M)}.
From Gauss' Divergence Theorem we have, 
($\nu $ denotes the exterior unit normal to $\partial \Omega $),
\begin{eqnarray*}
\label{gauss}
\int _{\Omega } x_N^\alpha|x|^{l} \, dx
 & = & 
\frac{1}{l+N+\alpha} \int_{\Omega } \mbox{div}\, \left( x_N^\alpha |x|^{l}x  \right) \, dx 
\\
 \nonumber
 & = & 
\frac{1}{l+N +\alpha} \int_{\partial \Omega } (x\cdot \nu ) x_N^\alpha |x|^{l} {\mathcal H}_{N-1}   ( dx) 
\\
 \nonumber
 & \leq & 
\frac{1}{l+N+\alpha } \int_{\partial \Omega } x_N^\alpha |x|^{l+1} {\mathcal H}_{N-1}   ( dx),
\end{eqnarray*}
with equality for $\Omega = B_{R}^{+} $,
and (\ref{ineqrad}) follows for smooth sets. 
Using (\ref{CklNsmooth}), we also obtain (\ref{ineqrad}) for measurable sets. 
\\
Finally assume $l+1<k$. Using Lemma \ref{rangekl1} and the result for $l+1=k$ 
we again obtain (\ref{ineqrad}), 
and (\ref{M=BR}) can hold only if $M= B_R^{+}$ for some $R>0$.
$\hfill \Box $ 

\noindent

\section{Application to a class of degenerate pde's}

In this section we show how our isoperimetric inequality allows to obtain  pointwise estimates for solutions to  a class of boundary value problems for 
degenerate elliptic
equations. We divide it into two parts: in the first one we prove a Polya-Szeg\"o principle and we deduce a  Poincar\'e-type inequality, while in the second part we provide the comparison result  given by Theorem \ref{Comp} below.

\subsection{P\'{o}lya-Szeg\"o principle}
First we obtain a P\'{o}lya-Szeg\"o principle related to our 
isoperimetric inequality (\ref{mainineq}) (see also \cite{Talenti2}) which applies
 to a special class of functions that we define below.


\noindent Assume that the numbers $N, k, l$ and $\alpha$ satisfy the assumptions of Theorem \ref{maintheorem}. Then (\ref{mainineq}) implies
\begin{equation} 
\label{Isop_klalpha}
\int_{\partial \Omega }
|x|^k x_N ^{\alpha } {\mathcal H}_{N-1}(dx)
\geq 
\int_{\partial \Omega ^{ \star }}
|x|^k x_N ^{\alpha } 
{\mathcal H}
_{N-1}(dx)
\end{equation}
for every smooth set $\Omega \subset \mathbb{R} ^N_+ $, where $\Omega ^{\star}$ is the $\mu_{l,\alpha } $-symmetrization of $\Omega $.
We will use (\ref{Isop_klalpha}) to prove the following result
\begin{theorem} 
(P\'{o}lya-Szeg\"o principle) 
Let the numbers $N, k,l$ and $\alpha $ satisfy the assumptions 
of Theorem \ref{maintheorem} and denote  $m:= 2k -  l $. Then there holds 
\begin{equation}
\int_{\Omega }
\left\vert  \nabla u\right\vert ^2 
d \mu _{m,\alpha } (x)
\geq 
\int_{\Omega^\star}\left\vert \nabla
u^{ \star }\right\vert ^{2}
d\mu_{m,\alpha } (x)\,,
\quad 
\forall u\in V^2(\Omega, d\mu_{m,\alpha}),
\label{PS_k_l}
\end{equation}
where $u^{\star } $ denotes the $\mu _{m,\alpha } $-symmetrization of $u$.
\end{theorem}
{\sl Proof:} 
It is sufficient to consider the case that $u$ is non-negative. Further,
by an approximation argument we may assume that
$u \in C^{\infty}_{0}(\mathbb{R}^{N} ) $.
Let 
\begin{eqnarray*}
I & := &  
\int_{\Omega } | \nabla u| ^{2}
|x| ^{2k-l} x_N ^{\alpha }\, dx \quad \mbox{and}\\
I ^{\star }  & := &  
\int_{\Omega^\star } | \nabla u^{\star} | ^{2}
|x| ^{2k-l} x_N ^{\alpha }\, dx .
\end{eqnarray*}
The coarea formula yields 
\begin{eqnarray*}
\label{1coarea}
I 
 & = & 
 \int_{0}^{\infty }\int_{u=t} |\nabla
u|  |x| ^{2k-l} x_N ^{\alpha }\, {\mathcal H}_{N-1}(dx)\, dt \quad \mbox{and}
\\
\label{coarea2}
I^{\star}
 & = & 
 \int_{0}^{\infty }\int_{u^{\star } =t} |\nabla
u^{\star} |  |x| ^{2k-l} x_N ^{\alpha }\, {\mathcal H}_{N-1}(dx)\, dt 
.
\end{eqnarray*}
Further, H\"older's inequality gives
\begin{equation}
\label{1holder}
\int_{ u=t} |x|^k x_N ^{\alpha } \, {\mathcal H} _{N-1} (dx)
\leq 
\left( \int_{ u=t} |x|^{2k - l} |\nabla u|  x_N ^{\alpha } \, {\mathcal H} _{N-1} (dx) \right) ^{\frac{1}{2} } 
\cdot
\left( \int_{ u=t} \frac{|x|^l x_N ^{\alpha }}{|\nabla u|} \, {\mathcal H}_{N-1} (dx) 
\right)
^{\frac{1}{2} } ,
\end{equation} 
for a.e. $t\in [0, +\infty )$. 
Hence (\ref{1coarea}) together with (\ref{1holder}) tells us that
\begin{equation}
\label{coarea3}
I 
\geq 
\int_{0}^{\infty }
\left( \int_{u=t} |x| ^{k} x_N ^{\alpha }\, {\mathcal H}
_{N-1}(dx)
\right) ^{2} \cdot 
\left( 
\int_{u=t}\frac{ |x| ^{l}x_N ^{\alpha }}{
| \nabla u| } x_N ^{\alpha } \, {\mathcal H}_{N-1}(dx)
\right) ^{-1} \, dt.
\end{equation}
Since $u^{\star} $ is a radial function, we obtain in an analogous manner,
\begin{equation}
\label{coarea4}
I^{\star} 
=
\int_{0}^{\infty }
\left( \int_{u^{\star} =t} |x| ^{k} x_N ^{\alpha } \, {\mathcal H}
_{N-1}(dx)
\right) ^{2} \cdot 
\left( 
\int_{u^{\star} =t}\frac{ |x| ^{l}x_N ^{\alpha } }{
| \nabla u^{\star} | } \, {\mathcal H}_{N-1}(dx)
\right) ^{-1} \, dt.
\end{equation}
Observing that
\begin{equation*}
\label{meas_u>t}
\int_{u>t} |x|^{l} x_N ^{\alpha } \, dx
=
\int_{u^{\star }>t}
|x|^{l} x_N ^{\alpha } \, dx \quad  \forall t\in [0, +\infty ), 
\end{equation*}
Fleming-Rishel's formula yields
\begin{equation}
\label{flemingrishel} 
\int_{u=t } \frac{|x|^l x_N ^{\alpha }}{|\nabla u|} \, {\mathcal H}_{N-1} (dx) 
=
\int_{u^{\star} =t } \frac{|x|^l x_N ^{\alpha }}{|\nabla u^{\star} |} \, {\mathcal H}_{N-1} (dx)
\end{equation}
for a.e. $t\in [0, +\infty )$. 
Hence
(\ref{flemingrishel}) and (\ref{Isop_klalpha}) give
\begin{eqnarray*}
 & &
\int_{0}^{\infty }
\left( \int_{u=t} |x|^k x_N ^{\alpha } \, {\mathcal H}
_{N-1}(dx) \right) ^{2}
\cdot 
\left( \int_{u=t}\frac{| x| ^{l} x_N ^{\alpha } }{
| \nabla u| } \, {\mathcal H}_{N-1}(dx) \right) ^{-1}
\, dt 
\\
 & \geq &
\int_{0}^{\infty }\left( \int_{u^{\star} =t} |x| ^{k} x_N ^{\alpha }
\, {\mathcal H}_{N-1}(dx) \right) ^{2} \cdot \left( \int_{u^{\star}=t}
\frac{|x|^{l} x_N ^{\alpha } }{|\nabla u^{\star} | } \, {\mathcal H}
_{N-1}(dx)\right) ^{-1} \, dt.
\end{eqnarray*}
Now (\ref{PS_k_l})  follows from this, (\ref{coarea3}) and (\ref{coarea4}).
$\hfill \Box$
\\[0.1cm]
%

\begin{corollary}
\label{poincare-dv} 
There exists a constant $C$, depending only on $\mu _{m,\alpha }(\Omega)$, such that 
for any $u\in V ^{2} 
(\Omega, d\mu _{m,\alpha }) $ the following inequality holds
\begin{equation*}
 \int_{\Omega}  u^{2}d\mu_{m,\alpha } \
 \leq C
\int_{\Omega}\left\vert \nabla u\right\vert ^{2}d\mu_{m,\alpha }  .  
\label{qp}
\end{equation*}
\end{corollary}

\noindent \textsl{Proof : } Let $u\in V ^{2} 
(\Omega , d\mu _{m,\alpha }) $ and define $U$ by $
U(|x|):=u^{\star }(x)$, ($x\in \Omega^{\star }$). Using (\ref{PS_k_l}), since rearrangements preserves the $L^{p}$-norms, we find 
\begin{eqnarray*}
\frac
{\displaystyle\int_{\Omega}\left\vert \nabla u\right\vert ^{2}d\mu_{m,\alpha } }
{\displaystyle\int_{\Omega}\left\vert u\right\vert ^{2}d\mu_{m,\alpha }  } 
&\geq &
\frac
{\displaystyle\int_{\Omega^\star}\left\vert \nabla u^{\star }\right\vert ^{2}d\mu_{m,\alpha } }
{\displaystyle\int_{\Omega^\star}\left\vert u^{\star }\right\vert ^{2}d\mu_{m,\alpha } }
 \\
&=&
\frac
{\displaystyle\int_{0}^{R^{\star }}\left\vert \frac{dU}{dr} \right\vert ^{2}r^{N-1+m+\alpha}\,dr}
{\displaystyle\int_{0}^{R^{\star }}U^{2}r^{N-1+m+\alpha}\,dr }
\geq 
\frac
{ \displaystyle\int_{0}^{R^{\star }}\left\vert \frac{dU}{dr}\right\vert ^{2}r^{N-1+m+\alpha}\,dr }
{\displaystyle\int_{0}^{R^{\star }}U^{2}r^{N-1+m+\alpha}\,dr}.
\end{eqnarray*}
Since $N+l+\alpha+2>0$ we can conclude the proof by applying Thorem 4 at p.45 in \cite{M2}.

\subsection{Comparison result}

Let us consider the class  of boundary value problems for degenerate elliptic
equations of the type
\begin{equation}
\left\{
\begin{array}{ccc}
-\text{div} (A(x)\nabla u)= & x_N^\alpha|x|^m f(x) & \text{in }\Omega \\
u=0\text{ \ \ \ \ \ \ \ \ \ \ \ } &  & \text{ \ on }\Gamma _{+} \\
\frac{\partial u}{\partial x_N } = 0 \text{ \ \ \ \ \ \ \ \ \ \ \ } &  &
\text{ \ on }\Gamma _{0}
\end{array}
\right.  \label{P}
\end{equation}
where $\Omega$ is an open subset of $\mathbb R^N_+$ with $\mu_{m,\alpha}(\Omega) < + \infty $, $\Gamma _+ =
\mathbb R^N_+ \cap \partial \Omega $, $\Gamma _{0} =  \partial \Omega \setminus \Gamma^+$, 
$m:=2k-l$. Assume that $A(x)=(a_{ij}(x))_{ij}$ is an $N\times N$ symmetric matrix with measurable
coefficients satisfying
\begin{equation}
x_N^\alpha|x|^m \left\vert \zeta \right\vert ^{2}\leq a_{ij}(x)\zeta _{i}\zeta
_{j}\leq \Lambda x_N^\alpha|x|^m \left\vert \zeta \right\vert ^{2},\text{ \ \ }
\Lambda \geq 1,  \label{ell}
\end{equation}
for almost every $x\in \Omega$ \ and for all $\zeta \in \mathbb{R}^{N}$. Assume
also that $f$ belongs to the weighted H\"older space $L^{2}(\Omega,d\mu_{m,\alpha} )$.

The type of degeneracy in (\ref{ell}) occurs, for $\alpha \in \mathbb{N}$, when
one looks for solutions to linear PDE's which are radially symmetric with
respect to groups of $(\alpha+1)$ variables (see, e.g.,  \cite{MadernaSalsa}
 and the references therein). The case of non-integers $\alpha$
has been the object of investigation, for instance, in the generalized
axially symmetric potential theory .

By a weak solution to such a problem \eqref{P} we mean a
function $u$ belonging to $V^2(\Omega ,d\mu_{m,\alpha} )$ such that
\begin{equation*}
\int_{\Omega }A(x)\nabla u\nabla \phi d\mu_{m,\alpha} =\int_{\Omega }f \phi d\mu_{m,\alpha} ,
\label{wsol-dv}
\end{equation*}
for every $\phi \in C^{1}(\bar{\Omega})$ such that $\phi =0$ on the set
$\partial \Omega \setminus \left\{ x_{N}=0\right\} $.

Observe that, by  Corollary  \ref{poincare-dv} , for any $f\in L^{2}(\Omega,d\mu_{m, \alpha})$, Lax-Milgram Theorem ensures the existence of a unique solution $u\in V ^{2}( \Omega  ,d\mu _{m, \alpha})$ to problem \eqref{P}.

We provide optimal bounds for the solution to problem
(\ref{P}) trough symmetrization methods introduced by G. Talenti in \cite{T}
(see also \cite{AT}, 
\cite{BBMP2} and
\cite{MadernaSalsa}).

\begin{theorem}
\label{Comp} Let $u$ be the weak solution to problem (\ref{P}), and let $w$
be the function
\begin{equation*}
w(x)=w^{\star }(x)=
\int_{\left\vert x\right\vert
}^{r^{\star }}\left( \int_{0}^{\rho }f^{\star }(\sigma )\sigma ^{N-1+m+\alpha
}d\sigma \right) \rho ^{-N+1-m-\alpha } d\rho ,
\end{equation*}
which is the weak solution to the problem
\begin{equation*}
\left\{ 
\begin{array}{cccc}
- \mathrm{div} \left( x_{N}^{\alpha }|x|^{m}\nabla w\right) = & x_{N}^{\alpha
}|x|^{m}f^{\star } & \mathrm{in } & \Omega ^{\star } 
\\ 
w=0\text{ \ \ \ \ \ \ \ \ \ \ \ } &  & \mathrm{on } & \partial \Omega ^{\star
}\setminus \left\{ x_{N}=0\right\}  \\ 
\frac{\partial w}{\partial x_{N}}=0  \text{ \ \ \ \ \ \ \ \ \ \ \ } &  &
 \mathrm{on } & \partial \Omega ^{\star }\cap \left\{ x_{N}=0\right\} 
\end{array}
\right. 
\label{Symm_Probl_dv}
\end{equation*}

Then
\begin{equation*}
u^{\star }(x)\leq w(x)\text{ a.e. in }\Omega^{\star },  \label{Point_Est-dv}
\end{equation*}
and
\begin{equation*}
\int_{\Omega}\left\vert \nabla u\right\vert ^{q}d\mu_{m, \alpha} \leq \int_{\Omega^{\star
}}\left\vert \nabla w\right\vert ^{q}d\mu_{m, \alpha} ,\text{ for all }0<q\leq 2.
\label{Grad-Est-dv}
\end{equation*}
\end{theorem}

The proof will be omitted since it can be easily obtained by repeating the
arguments of analogous results contained for instance in \cite{AT}, 
\cite{BBMP2}, 
\cite{BCM2} and \cite{T}.


\end{document}